%% file: flat.tex
\newcommand{\N}{\mathbb{N}}               

\newcommand{\R}{\mathbb{R}}

\newcommand{\A}{\mathscr{A}}

\newcommand{\m}{\mathfrak{M}}  
   
\newcommand{\E}{{\mathscr E}}   
\mathchardef\varepsilon="010F
\mathchardef\epsilon="0122
\mathchardef\vartheta="0112
\mathchardef\theta="0123
\mathchardef\varrho="011A
\mathchardef\rho="0125
\mathchardef\varphi="011E     
\mathchardef\phi="0127
\renewcommand \emptyset \varnothing
\documentclass[12pt]{article}            
\usepackage[all]{xy}
\usepackage[latin1]{inputenc}        
\usepackage[dvips]{graphics}
\usepackage[dvips]{graphicx}
\usepackage{amsfonts}
\usepackage{amssymb}
\usepackage{amsmath}
\usepackage{amstext}
\usepackage{amsbsy}
\usepackage{amsopn}
\usepackage{amsthm}
\usepackage{amscd}
\usepackage{amsxtra}
\usepackage{mathrsfs}
\usepackage{upref}
\usepackage{fullpage} 

\author{Gianmarco Capitanio\footnote{Research supported by INdAM.} 
}
\title{Cusp singularities of plane envelopes}

\begin{document}        

\theoremstyle{plain}
\newtheorem*{theorem*}{\bf Theorem}
\newtheorem{theorem}{\bf Theorem}
\newtheorem{lemma}{\bf Lemma}
\newtheorem*{lemma*}{\bf Lemma}
\newtheorem{proposition}{\bf Proposition}
\newtheorem*{proposition*}{\bf Proposition}
\newtheorem{corollary}{\bf Corollary}
\newtheorem*{corollary*}{\bf Corollary}
\theoremstyle{definition}
\newtheorem*{definition*}{\bf Definition}
\newtheorem*{definitions*}{\bf Definitions}
\newtheorem*{conjecture}{\bf Conjecture}
\newtheorem{example}{\bf Example}
\newtheorem*{example*}{\bf Example}
\theoremstyle{remark}
\newtheorem*{remark*}{\bf Remark}
\newtheorem{remark}{\bf Remark}
\newtheorem*{acknowledgements}{\bf Acknowledgements}
\maketitle
\begin{abstract}
We consider smooth $1$-parameter families of plane curves tangent to a
semicubic parabola, when the curvature radius of their curves at the
tangency point vanishes at the cusp point. 
We find the $\A$-normal form of these families, their envelopes and
local patterns near the cusp. 
We obtain a new codimension $2$ singularity of envelopes (two
transversal semicubic cusps), and  we describe its perestroikas under
generic small deformations.  
\end{abstract}
\section{Introduction}

In \cite{problems}, problem 2001-6, V. I. Arnold proposed to study envelopes of families of curves tangent to a semicubic parabola, when the curvature radius at the tangency point is a non-negative smooth function of this point, vanishing at the cusp point. 
This problem has been solved in \cite{cusp}, where such an envelope is described, according to the value of the first term of the Taylor expansion of the curvature radius function. 
These families are not smooth. 

In this paper we consider a similar problem in the smooth setting. 
Namely, we consider {\it cusped tangential families}, that is smooth $1$-parameter families of smooth plane curves, tangent to a semicubic parabola  (called the support). 
This is a generalization of the notion of tangential family (see \cite{stable}, \cite{simple} and \cite{dual}) to the case of singular support.

If the curvature radius of the family curve tangent to the support at the singular point is not vanishing, then the curve is regular.  
Such families provide one of the five local models for generic $1$-parameter families of plane curves (studied by Dufour, see e.g. \cite{dufour1}, \cite{dufour2}).   

The aim of this paper is to study cusped tangential families when this
curvature radius vanishes (we say that such a family is flat along the
cusp). 
This means that we allow the family to have a singular curve, corresponding to the singular point of the support. 
These families provide a new codimension $2$ local singularity of envelopes.  
While stable singularities of plane envelopes have been fully classified (by Arnold \cite{arnold1976b}, Dufour  \cite{dufour1} and Thom \cite{thom}), simple singularities have not yet been completely investigated, in spite of the fact that they appear generically in deformations of families with enough parameters.  

We describe the singularities of envelopes of flat cusped tangential
families and their local patterns: we show that, generically, 
the envelope has two branches near the cusp of the support, the support itself and a second semicubic cusp transversal to it. 
This singularity is local, hence it is different from the superposition of two branches having both a semicubic cusp. 
Near the support singular point, the curves of the family experience a $\gamma\rightarrow U$ perestroika along the support.      
We discuss the graph surfaces of cusped tangential families (whose apparent contour is the family envelope). 
Moreover, we find $\A$-normal forms of generic cusped tangential families, toghether with their $\A$-miniversal deformations. 
This allows us to describe the perestroikas of the envelope singularity under small deformations. 
These perestroikas can be interpreted as the metamorphosis of the
apparent contour of a vertical Whitney Umbrella under small
deformations of the direction of the projection. 
 
\section{Cusped tangential families} \label{sct:2}

Let $\Gamma$ be a smooth curve in the plane $\R^2=\{X,Y\}$, having a semicubic cusp. 
Without loss of generality, we may assume that, at least near the
singular point, this curve is 
the standard semicubic parabola $Y^2=X^3$, parameterized by the map 
$\xi\mapsto \Gamma(\xi)=(\xi^2,\xi^3)$.    

\begin{definition*}
A smooth $1$-parameter family of plane curves $\{\phi_\xi:\xi\in\R\}$ is a {\it cusped tangential family} (CTF) of support $\Gamma$ if each curve     $\phi_\xi:\R\rightarrow\R^2$ is tangent to the semicubic parabola $\Gamma$ at $\Gamma(\xi)$ for every $\xi\not=0$.   
\end{definition*}

By ``smooth family'' we mean that the mapping $\phi:\R\times \R\rightarrow\R^2$, sending $(\xi,t)$ to $\phi(\xi,t):=\phi_\xi(t)$, is of class $C^\infty$.  
In particular, each curve $\phi_\xi$ is smooth. 

Since the semicubic parabola $\Gamma$ is regular at each point $\xi\not=0$, each curve on the family, corresponding to such a point, is required to be regular near the tangency point. 
Up to change the parameterizations of the tangent curves, we can suppose that they are tangent to $\Gamma$ for $t=0$ at point $\phi_\xi(0)=\Gamma(\xi)$.  

On the other hand, the special curve $\phi_0$, corresponding to the singular point of the semicubic parabola, may have a singularity at $t=0$. 

The vector field along the semicubic parabola $\Gamma$, 
\begin{equation}\label{eq:1}
V_\Gamma(\xi):=\dfrac{1}{\xi}\dfrac{d}{d\xi}\Gamma(\xi)
\end{equation}
is smooth also at the singular point $\xi=0$, and it defines a tangent direction at each point of the semicubic parabola. 
The tangency condition between $\phi_\xi$ and $\Gamma$ can be expressed as the proportionality of the tangent vector field  $V_\Gamma$ and the vector field of the velocities of the curves $\phi_\xi$ at the tangency point:
\begin{equation}\label{eq:2}
\dfrac{d\phi_\xi}{dt}(0) = \alpha(\xi) \ V_\Gamma(\xi)\ ,  
\end{equation}
where the {\it proportionality factor} $\alpha:\R\rightarrow\R$ is a smooth function, nonvanishing at every $\xi\not=0$.  
This proportionality factor $\alpha$ does not characterize a CTF, but it controls several interesting features of the family. 

The degeneracy of the proportionality factor at the singular point of
the support provides a stratification of the set of the CTFs.

\begin{definition*} 
A CTF is said to be {\it $n$-flat} along $\Gamma$ if its proportionality factor $\alpha$ is $n$-flat at the origin (i.e. the Taylor expansion of $\alpha$ at     $\xi=0$ is $O(\xi^n)$).
\end{definition*} 

$0$-flat CTFs will be also called {\it non-flat} CTFs. 
By {\it flat} CTF we mean a CTF whose proportionality factor vanishes at $\xi=0$. 
There exist infinitely degenerate flat families which are not
$n$-flat, whatever be $n\in\N$.  

\begin{remark*}
A CTF is non-flat if and only if the special curve $\phi_0$ of the family is regular at $t=0$ and it is tangent there to the support. 
For flat CTFs, the special curve is singular at the singular point of
the semicubic parabola $\Gamma$.  
In particular, its curvature radius is vanishing at the cusp point of
the support.
\end{remark*} 

We introduce now a ``genericity condition'' on CTFs. 

\begin{definition*} 
We say that a flat CTF satisfies the {\it condition $(*)$} if its special curve $\phi_0$ has a semicubic cusp at $t=0$, transversal to that of   $\Gamma$ at the common singular point. 
\end{definition*} 

We will explain below (corollary to Theorem \ref{thm:3} and subsequent remark) why condition $(*)$ is actually a genericity condition.  

\section{Normal forms of $1$-flat CTFs}
In this section we find the $\A$-miniversal normal form of $1$-flat $(*)$-generic CTFs and its $\A$-miniversal deformation.
We start with some preliminar computations. 

Let $\phi$ be a parameterization of a CTF in the standard coordinate
system described above.   
Equation \eqref{eq:2} implies that the Taylor expansion of $\phi_\xi$ (as a function of the only variable $t$) at $t=0$ can be written as  
\begin{equation}\label{eq:3}
\phi_\xi(t)= \Gamma(\xi) + \alpha(\xi) V_\Gamma(\xi) t + 
\begin{pmatrix}
A(\xi) \\B(\xi)
\end{pmatrix}
t^2 + 
\begin{pmatrix}
C(\xi) \\D(\xi)
\end{pmatrix}
t^3 + 
\begin{pmatrix}
1 \\1
\end{pmatrix}
o(t^3) \ ,
\end{equation}
where $A,B,C,D:\R\rightarrow\R$ are smooth functions of the parameter $\xi$. 
We shall use the notation $\sum_{i=0}^\infty A_i\xi^i$ for the expansion of the function $A$ at $\xi=0$ (and similar notations for functions $B,C$ and $D$). 

\begin{lemma}\label{lemma:1}
A flat CTF $\phi$ is $(*)$-generic if and only if the coefficients of the Taylor expansion \eqref{eq:3} satisfy $B_0(A_0D_0-B_0C_0)\not=0$. 
\end{lemma} 

\begin{proof}
In our particular coordinate system, the $3$-jet at the origin of the special curve of the family is
$(A_0t^2+C_0 t^3, B_0t^2+D_0 t^3)$, since $\alpha(0)=0$. 
This curve has a semicubic cusp if and only if $A_0D_0\not=B_0C_0$. 
Its tangent direction at the singular point, defined by the tangent vector field \eqref{eq:1}, is $(A_0,B_0)$. 
It is transverse to that of the semicubic parabola $\Gamma$ if and only if $B_0\not=0$.  
\end{proof}

We shall classify $(*)$-generic $1$-flat CTFs with respect to the usual $\A$-equivalence. 
Let us recall that two smooth map germs $f,g:(\R^2,0)\rightarrow(\R^2,0)$ are said to be $\A$-equivalent if there exist diffeomorphic coordinate changes
$h,k:(\R^2,0)\rightarrow(\R^2,0)$ in the source and target space 
such that $g=k\circ f\circ h^{-1}$. 
In other terms, two germs are $\A$-equivalent if they belong to the same orbit of the natural action of the group $\A=\{(h,k)\}$ on  $\m_2\E_2^2$, where $\m_2$ is the maximal ideal of ring $\E_2$ of the smooth function germs $(\R^2,0)\rightarrow \R$ and $\E_2^2=\E_2\times\E_2$.

\begin{remark*}
$\A$-equivalence does not preserve the curves of the family, but only their envelopes. 
The more natural equivalence relation, preserving the family curves, is the fibered equivalence, restricting the source coordinate changes to those of the form  
$h(\xi,t)=(h_1(\xi),h_2(\xi,t))$.  
We will use coordinate $x,y$ in the source space, instead of $\xi,t$,
when the fibered structure of the source space is not taken into
account. 

Non-flat CTFs are local models (for this fibered equivalence) 
for semicubic cusps as stable singularities of envelopes of
$1$-parameter families of plane curves, studied by J.-P. Dufour (see
e.g. \cite{dufour1}, \cite{dufour2} and \cite{dufour3}).  
Among other results, he proved that 
any non-flat CTFs has a functional modulus, which is 
essentially the {\it Blaschke curvature} of
the (codimension $1$) planar $3$-web defined in the interior of the
semicubic parabola near the singular point (where at each point passes
$3$ curves of the family). 
Since flat CTFs are degenerated non-flat CTFs, it is clear that they
have fonctional moduli.  
Hence the fibered equivalence would not provide a discrete
classification of CTFs.  
\end{remark*} 

Some easy computations, starting from parameterization \eqref{eq:3}, show that the $3$-jet of $\phi$ is
$\A$-equivalent to 
$$\psi_\delta(x,y):=(x^2+y^2+\delta y^3,y^2+x^3) \ , $$ 
where $\delta$ is a parameter determined (up to sign) by the initial CTF. 

\begin{theorem}\label{thm:1}
Every $(*)$-generic $1$-flat CTF is $\A$-equivalent to the normal form
$\psi_\delta$, provided that $\delta\not=0$. 
\end{theorem}

For the proof of this theorem we need to recall some further
notations. 
The tangent space $V(f)$ to $\m_2\E_2^2$ at any one of its points $f$ is canonically identified  with $\m_2\E_2^2$ itself, and it contains the tangent space $T\A(f)$ to the $\A$-orbit of $f$. 
We recall that the tangent space is 
$$T\A(f):=\m_2 \langle f_{x},f_{y} \rangle_\R + f^*(\m_2) \langle
e_1,e_2\rangle_\R \ ,$$ 
where $f_{x}$ (resp. $f_y$) is the partial derivative of $f$ with
respect to $x$ (resp. with respect to $y$) and $e_1=(1,0)$, $e_2=(0,1)$. 

The {\it extended tangent space} to the $\A$-orbit of $f$, $T_e\A(f)$,
is the extension of $T\A(f)$ obtained by including the initial
velocity vectors of all $1$-parameter deformations of mappings and of
the identity diffeomorphisms  of the source and the target spaces:  
$$T_e\A(f):=\E_2 \langle f_{x},f_{y}  \rangle_\R + f^*(\E_2) \langle
e_1,e_2\rangle_\R \ ,$$   
The {\it extended $\A$-codimension} of $f$ is the dimension of the
real vector space $\E_n^p/T_e\A(f)$.  
We denote by $I_f$ the ideal $\E_2^2 \langle f_1,f_2\rangle_\R$
generated by the components of $f=(f_1,f_2)$.
 
The proof of Theorem \ref{thm:1} is based on  the following
determinacy estimate (see \cite{duplessis}).

\begin{theorem*}[Du Plessis] 
Let $f:(\R^2,0)\rightarrow  (\R^2,0)$ be a smooth map germ and suppose that 
\begin{equation}\label{eq:4}
\m_2^\ell \E_2^2 \subseteq \E^2 \langle f_x,f_y\rangle_\R + 
I_f\E_2^2 + \m_2^{\ell+1}\E_2^2
\end{equation}
and 
\begin{equation}\label{eq:5}
\m_2^k \E_2^2 \subseteq T_e\E(f) + \m_2^{k+\ell}\E_2^2
\end{equation}
for some integers $k\geq 1$ and $\ell>0$. 
Then $f$ is $(k+\ell)$-$\A$-determined.
\end{theorem*}

\begin{proof}[Proof of Theorem \ref{thm:1}]
We first use du Plessis' estimate with $k=\ell=2$ to show that the
normal form $\psi=(\psi_1,\psi_2)$ is $4$-$\A$-determined, provided
that $\delta\not=0$ (we omit here the subscript $\delta$).

Consider the $8$ vectors $x\psi_x,\dots,y \psi_y$, $(\psi_1,0)$,
$(\psi_2,0)$, $(0,\psi_1)$ and $(0,\psi_2)$ of the space 
$\E^2 \langle \psi_x,\psi_y\rangle_\R + I_\psi\E_2^2$. 
Their projections into the $6$-dimensional vector space
$\m_2^2\E_2^2/\m_2^3\E_2^2$, generated by the vector monomials
$(x^2,0),\dots,(0,y^2)$, form a generator system. 
Indeed, the matrix formed by the coordinates of these projections in
the above space basis has maximal rank. 
Hence 
$\E^2 \langle \psi_x, \psi_y \rangle_\R + 
I_\psi\E_2^2 + \m_2^{3}\E_2^2$ contains $\m_2^2 \E_2^2$, so the first
inclusion \eqref{eq:4} is fullfilled. 

The second inclusion \eqref{eq:5} to check is: 
$$\m_2^2 \E_2^2 \subseteq T_e\A(\psi) + \m_2^{4}\E_2^2\ .$$
Consider the $6$ vectors $x^2\psi_x,\dots,y^2\psi_y$, toghether with
the $8$ vectors used in the proof of the first inclusion.   
These $14$ vectors belong to the extended tangent space $T_e\A(\psi)$. 
Consider the squared matrix, whose columns are the coordinates of
the projections of these $14$ vectors in the space
$\m_2^2\E_2^2/\m_2^4\E_2^2$, equipped with the standard basis 
$$\{(x^i y^j,0),(0,x^i y^j) : i+j=2,3\} \ .$$
The determinant of this matrix is $1280\ \delta$, so 
the above inclusion holds whenever $\delta\not=0$.

By Du Plessis' estimate, $\psi$ is $4$-$\A$-determined. 
Now, all the degree $4$ terms of $\psi$ can be killed without any change of the $3$-jet of the germ (provided that $\delta\not=0$), so the germ $\psi$ is $3$-$\A$-determined. 
\end{proof} 

\section{Envelopes of flat CTFs}\label{sct:4}
In this section we describe envelopes of flat CTFs.  
First let us recall some basic definitions of envelope theory. 

The {\it graph} of a family of curves  $\{\phi_\xi:\xi\in\R\}$ is the surface 
$$\left\{\big(\phi_\xi(t),\xi\big) : \xi,t\in\R\right\}$$
in the $3$-space $\R^2\times\R=\{X,Y,Z\}$.
The {\it envelope} of the family is the apparent contour of its graph on the plane $\R^2$ (i.e., the critical value set of the restriction to the graph of the natural projection ``forgetting $Z$'' 
$\R^2\times\R\rightarrow\R^2$). 
There exist several other definitions of envelope, but the one we recalled (due to René Thom \cite{thom}) is the most general and  geometric. 

\begin{remark*}
By the very definition of envelope,  the semicubic support of every CTF  is a branch of its envelope. 
In the case of germs of non-flat CTFs, it is the only branch of the envelope. 
\end{remark*}

For germs of flat CTFs the envelope is more complicated. 

\begin{theorem}\label{thm:2}
The envelope of the germ at the origin of any $(*)$-generic flat CTF has, in addition to the semicubic support $\Gamma$, just another branch, having a cusp at the origin.  

For $1$-flat CTFs, this second cusp is semicubic if and only if the parameter $\delta$ of its $\A$-normal form is not vanishing.
Moreover, this cusp is transversal to those of the support $\Gamma$ and of the special curve $\phi_0$.

For $(n>1)$-flat CTFs, the second envelope cusp is semicubic and tangent to that of the special curve $\phi_0$.
\end{theorem}

\begin{proof}
Consider first the case of $1$-flat CTFs. 
Let $\phi:\R^2\rightarrow\R^2$ be such a CTF. 
By Theorem \ref{thm:1}, $\phi$ is $\A$-equivalent to the normal form $\psi_\delta$. 
Hence, the critical value sets of $\phi$ and $\psi_\delta$ are diffeomorphic. 
The Jacobian determinant of $\psi_\delta$ is 
$4xy-3(2+3\delta y)x^2y$.
Therefore the critical set has two smooth branches passing through the origin, of equations $x=0$ and $y=0$. 
The two branches of the critical value sets are $x\mapsto(x^2,x^3)$,
which is the support $\Gamma$, and $y\mapsto (-y^2+\delta y^3,y^2)$. 
This second cusp is semicubic if and only if $\delta\not=0$.

Consider now an $n$-flat CTF, for $n\geq 2$. 
A direct computation shows that the $2$-jet of the Jacobian determinant of the parameterization \eqref{eq:3} is $4B_0\xi t$.  
Since $B_0$ is non-zero by Lemma \ref{lemma:1}, the critical point equation has, in addition to the solution $t=0$ 
(providing the envelope branch $\Gamma$), a solution of the form $\xi=o(t)$. 
Replacing in \eqref{eq:3}, we find the envelope parameterization's $3$-jet to be equal to that of the special curve $\phi_0$.    
\end{proof}

\section{Singularities of graphs of flat CTFs}\label{sct:5}

In order to understand the behaviour of the curves of flat CTFs near the singular points of their supports, it is useful to describe the singularities of their  graphs. 
For this, we identify such a graph to the image of the mapping 
$$\Phi: \R^2\rightarrow\R^3\ , \quad  (\xi,t)\mapsto \Phi(\xi,t)= \big(\phi(\xi,t),\xi\big) \ . $$
We are interested to the classification of graphs under
$\A$-equivalence of map germs $(\R^2,0)\rightarrow (\R^3,0)$ ---the
definition of $\A$-equivalence in this case is similar to that
recalled in section \ref{sct:2}. 

Local singularities of mappings from the plane to the $3$-space have been classified under $\A$-equivalence by Whitney (who proved that the only stable singularities of these mappings are  transversal intersections of two or three regular sheets and the so-called Whitney Umbrella, see \cite{whitney}) and Mond (who classified simple singularities of these mappings, see \cite{mond}). 

Let us recall that a smooth map from the plane to the space has an $A^\pm_n$ singularity at, say, the origin if it is locally $\A$-equivalent to the normal form 
\begin{equation}\label{eq:An}
(u,v)\mapsto (u,v^2,v^3\pm u^n v)\ .
\end{equation} 
$A_n^+$ and $A_n^-$ singularities are different if and only if $n$ is even. 

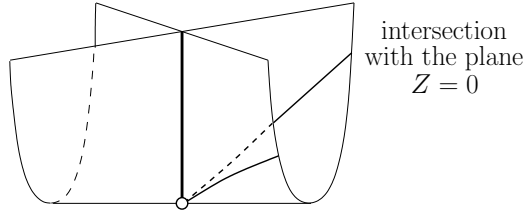
\begin{figure}[h]
  \centering
   \scalebox{.6}{\input{fig_1_flat.pstex_t}}
  \caption{The image of the Whitney Umbrella singularity and its
   intersection with a transversal plane.}
 \label{fig:1}
\end{figure}

$A_1$ singularity is stable; it is also called the {\it Whitney Umbrella} singularity ---due to the form of its image, see fig. \ref{fig:1}.    
The intersection of the (image) surface of the $A^\pm_n$-singularity with a transversal plane passing through the singular point is a semicubic cusp. 
$A_n^\pm$ singularities are simple for every $n\in\N$ (that is, by deforming them slightly we can get only a finite number of different singularities). 

\begin{theorem}\label{thm:3}
The graph of the germ of any $(*)$-generic $n$-flat CTF has a singularity $A^\pm_n$. 
\end{theorem}

\begin{proof}
By equation \eqref{eq:3}, the $(n+1)$-jet of any $n$-flat CTF $\phi$ can be
written, in a suitable coordinate system, as 
$$j^{n+1}\phi(\xi,t) =\Big( \xi^2+2\alpha_n\xi^n t+\sum
P_{i,j}\xi^i t^{j+2} , \xi^3+ \sum  Q_{i,j}\xi^i
t^{j+2},\xi  \Big) \ , $$
where both sums are taken over all the indexes $i,j\in\N\cup\{0\}$ such that
$i+j\leq n-1$. 
According to our preceding notations \eqref{eq:3}, $P_{0,0}=A_0$, $P_{0,1}=C_0$, $Q_{0,0}=B_0$, $Q_{0,1}=D_0$ and so on. 
By $(*)$-genericity, the coefficient $Q_{0,0}=B_0$ is not vanishing due to Lemma \ref{lemma:1}.  
Hence, it is easy to bring the above jet, by an $\A$-conjugation, to the form
$$\Big( 
2\alpha_nx^n y+\sum_{i+j\leq n-1} P_{i,j}x^i y^{j+2} ,
y^2  + \sum_{1\leq i+j\leq n-1} \tilde Q_{i,j}x^i y^{j+2},x 
\Big) \ , $$
where $\tilde Q_{i,j}=Q_{i,j}/Q_{0,0}$. 
The coordinate change $(X,Y,Z)\mapsto (X-P_{0,0}Y,Y,Z)$ in the target space of the map reduces the latter jet to 
$$\Big( 
2\alpha_nx^n y+\sum_{1\leq i+j\leq n-1} \tilde P_{i,j}x^i y^{j+2} ,
y^2  + \sum_{1\leq i+j\leq n-1} \tilde Q_{i,j}x^i y^{j+2},x 
\Big)\ , $$
where $\tilde P_{i,j}=P_{i,j}-P_{0,0}\tilde Q_{i,j}$. 
In particular, $\tilde P_{0,1}=P_{0,1}-P_{0,0}Q_{0,1}/Q_{0,0}$ 
vanishes if and only if $P_{0,1} Q_{0,0}=P_{0,0}Q_{0,1}$, i.e. if and only if  $B_0 C_0=A_0 D_0$, which is excluded by Lemma \ref{lemma:1}.

Now, by the coordinate change 
$$(X,Y,Z)\mapsto (X,Y+\sum_{k\geq 0} \mu_k\ Y Z^k,Z)$$ 
in the target space, we can kill (for suitable values $\mu_k$) all the
monomials $\{x^i t^{j+2} : i\geq 1\}$ in the $Y$-coordinate of the germ. 
Hence the initial jet is $\A$-equivalent to 
\begin{align*}
\Big(2\alpha_nx^n y+&\sum_{i=1}^{n-1} \tilde P_{i,0}x^i y^2 + \tilde P_{0,1} \ y^3
+ \sum_{i=1}^{n-2} \tilde P_{i,1}x^i y^3 + 
\sum_{i+j\leq n-3} \tilde P_{i,j}x^i y^{j+4}, \\
&y^2  + \sum_{i+j\leq n-2} \tilde Q'_{i,j}x^i y^{j+3},x  \Big) \ . 
\end{align*}
Acting similarly on the first coordinate of the germ, we reduce it to: 
$$\Big(  
2\alpha_nx^n y + \tilde P_{0,1} \ y^3
+ \sum_{i=1}^{n-1} \tilde P'_{i,2}x^i y^3 + 
\sum_{i+j\leq n-3} \tilde P'_{i,j}x^i y^{j+4}, 
t^2  + \sum_{i+j\leq n-2} \tilde Q'_{i,j}x^i y^{j+3},x 
\Big) \ . $$

For a suitable choice of parameters $\nu_{i,j}$, replacing $y$ with 
$y+\sum \nu_{i,j} x^i y^j$, $i+j\geq 2$, we get  
$$\Big(  
2\alpha_nx^n y + \tilde P_{0,1} \ y^3 + \sum_{i=1}^{n-2} \hat P_{i,1}x^i y^3 +  \sum_{i+j\leq n-3} \hat P_{i,j}x^i y^{j+4}, 
y^2, x 
\Big) \ . $$
Since $\tilde P_{0,1}$ is not vanishing, by rescaling we get 
$$\left(  
y^3 \pm x^n y + \sum_{i=1}^{n-2} \bar P_{i,1}x^i y^3 +  \sum_{i+j\leq n-3} \bar P_{i,j}x^i y^{j+4}, 
y^2, x  \right) \ , $$
which is clearly $\A$-equivalent to the normal form of the $A^\pm_n$ singularity \eqref{eq:An}; the sign $\pm$ is the same sign as the product $\alpha_n(B_0 C_0-A_0D_0)$ of the coefficients of the initial map \eqref{eq:3}.
This jet is $(n+1)$-$\A$-sufficient (see \cite{mond}), so the initial
CTF's graph is $\A$-equivalent to it. 
\end{proof}

We shall focalize now on $1$-flat CTFs. 
In this case Theorem \ref{thm:3} provides the following. 

\begin{corollary*}
The graph of the germ of any $(*)$-generic $1$-flat CTF has a Whitney Umbrella
singularity, transversal to the plane $Z=0$.
\end{corollary*}

\begin{remark*}
This fact explains why we consider condition $(*)$ a genericity condition: it expresses the genericity of the intersection of the graph of the family with the plane $Z=0$ passing through the singular point. 
\end{remark*}
 
\begin{remark*}
It follows from the corollary that the curves of any $(*)$-generic $(n=1)$-flat CTF experience a $\gamma\rightarrow U$ perestroika along the support  
$\Gamma$.
In the case of more flat CTFs, we can still deduce the pattern of the curves of the family near the singular point of the semicubic support. 
For $n=2$, we have two different cases, according to the sign $\pm$ of the singularity $A^\pm_2$ of the graph of the family. 
For the minus sign, the curves experiences a ``$\gamma\rightarrow\gamma$ perestroika along the cusp $\Gamma$", while for the plus sign, they experiences a ``$U\rightarrow U$ perestroika along the cusp $\Gamma$". 
For $n>2$ the patterns are similar to the preceding ones, according to the parity of the flatness degree.
\end{remark*}

\begin{figure}[h]
  \centering
   \scalebox{.31}{\input{fig_2_flat.pstex_t}}
  \caption{Local patterns of $(*)$-generic $1$-flat CTFs.}
 \label{fig:2}
\end{figure}
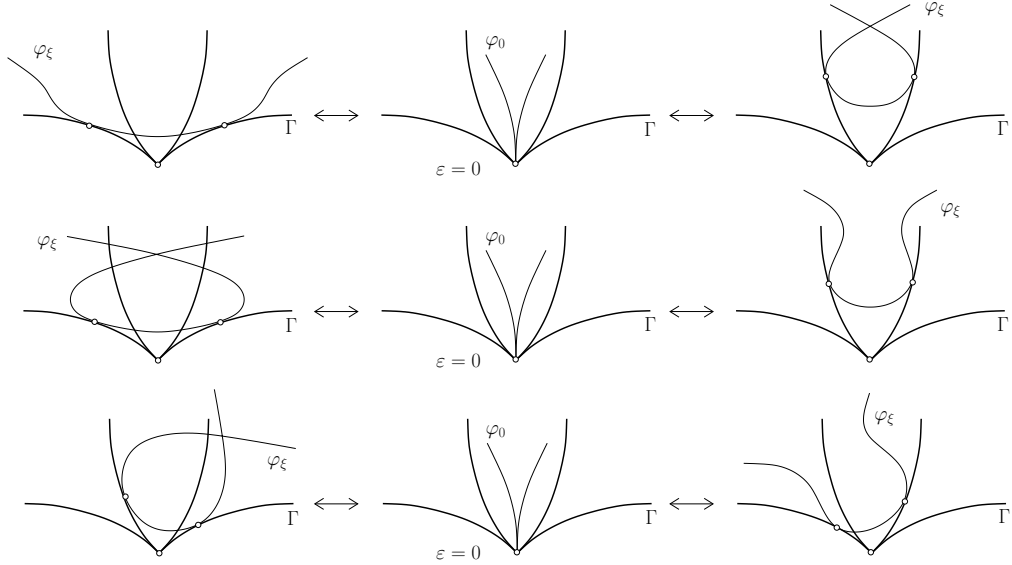

The local patterns of $(*)$-generic $1$-flat CTFs near the double cusp
singularity of their envelopes are shown in figure \ref{fig:2}.

\section{Deformations of the double cusps singularity}
We descibe now the metamorphosis occurring to the double cusp envelope singularity under a small deformation. 
This is done by computing the critical value sets of an $\A$-miniversal deformation of the normal form $\psi_\delta$.

\begin{theorem}\label{thm:4}
The map germ $\Psi_\delta:(\R^2\times\R^3,0)\rightarrow(\R^2,0)$ defined by 
$$ \Psi_\delta (x,y;\lambda,\mu,\nu):= \psi_\delta (x,y)+(\lambda y+\nu y^3,\mu x) $$
is an $\A$-miniversal deformation of the normal form $\psi_\delta$,
provided that $\delta\not=0$.
\end{theorem}

\begin{proof}
Recall that for every $r$-$\A$-determined map germ
$f:(\R^2,0)\rightarrow(\R^2,0)$ one has 
$$\m_2^{r+1} \E_2^2\subseteq T\A(f)\subseteq T_e\A(f)$$
(see \cite{gaffney}). 
Hence, the real vector space $\E_2^2/T_e\A(\psi_\delta)$ is generated by the vector monomials $(x^i y^j,0)$, $(0,x^i y^j)$, for all the non-negative integers $i,j$ such that $i+j=1,2,3$.  
Some easy computations lead to the equality  
$$\E_2^2=T_e\A(\psi_\delta) \oplus \R\left\{ 
\left(\begin{smallmatrix} y \\ 0 \end{smallmatrix}\right) , 
\left(\begin{smallmatrix} 0 \\ x \end{smallmatrix}\right) , 
\left(\begin{smallmatrix} y^3 \\ 0 \end{smallmatrix}\right) 
\right\} \ , $$
proving the theorem.
\end{proof}

\begin{corollary*}
The extended codimension of any $(*)$-generic $1$-flat CTF, whose coefficient $\delta$ is not vanishing, is  $3$.
\end{corollary*}

\begin{remark*}
The double cusp envelope singularity of $(*)$-generic $1$-flat CTFs is a codimension $2$ envelope singularity. 
Indeed the parameter $\nu$ in the miniversal deformation $\Psi_\delta$ does not change the envelope singularity, provided that $\delta\not=0$ and $\nu$ is small enough.
Hence, this singularity appear generically in envelopes of
$1$-parameter family of plane curves depending in two external
parameters, and it is not avoidable with arbitrarily small
deformations. 
\end{remark*}

We shall consider  from now on only the $2$-parameter deformation $\tilde \Psi_\delta:=\Psi_\delta|_{\nu=0}$ of the normal form $\psi_\delta$. 
We start considering the two sub-deformations 
$$H(x,y;\lambda):= \tilde \Psi_\delta(x,y;\lambda,0)  \ \text{and} \ K(x,y;\mu):= \tilde \Psi_\delta(x,y;0,\mu) \ ,$$ 
obtained by fixing $\mu=0$  and $\lambda=0$ respectively.

\begin{theorem}\label{thm:5}
Under deformations $H$ and $K$, one of the two cusps of the envelope singularity is preserved (being diffeomorphic to the standard semicubic parabola for every small enough value of the deformation parameter), while the other cusp experiences a $\gamma\rightarrow U$ transition when the deformation parameter crosses $0$. 
\end{theorem}

\begin{remark*}
The semicubic cusp is a stable envelope singularity: an isolated
envelope cusp cannot experience a $\gamma\rightarrow U$ transition. 
\end{remark*}

\begin{proof}
We prove the theorem for deformation $H$. 
The proof for $K$ is similar (the role of the two cusps being exchanged).

The vanishing of the Jacobian determinant (with respect to $x,y$) of $H$ provides the critical set equation
\begin{equation}\label{eq:7}
3x^2(\lambda+2y+3\delta y^2)=4xy \ ,
\end{equation} 
having two solutions passing through the origin $(0,0)$. 
The first solution $x=0$ provides a $\lambda$-depending family of critical value set branches.
The graph of this family in the $3$-space $\R^3=\{X,Y,Z\}$ is a surface, parameterized near the origin  by 
$(\lambda y +y^2+\delta y^3,y^2,\lambda)$. 
Since we are assuming $\delta\not=0$, this map germ is clearly $\A$-equivalent to the Whitney Umbrella normal form $(\lambda y,y^2,\lambda)$. 
Moreover, the intersection of this surface with the plane $Z=0$ is transversal, because the critical value set branch has a semicubic cusp at the origin for $\lambda=0$. 
Hence these critical value set branches have a $\gamma\rightarrow U$ transition when $\lambda$ crosses $0$. 

The other relevant solution of \eqref{eq:7} is also a solution of  
\begin{equation}\label{eq:8}
(9\delta x)y^2+(6x-4)y+3\lambda x \ .
\end{equation} 
Let $\Delta$ denote the discriminant of this second degree equation in
the unknown $y$, as a function of the variables $x,\lambda$ ($\delta$
being fixed); since $\Delta(x=0,\lambda=0)=4$, $\sqrt\Delta$ is a smooth function germ (at $x=\lambda=0$).
A straght-forward computation shows that 
$$\sqrt{\Delta}= 4-6x-\frac{27}{2} \delta \lambda x^2+ o_{\lambda,x}(3) \ , $$
where $o_{\lambda,x}(n)$ is some smooth function germ of order higher
than $n$ in the variables $\lambda,x$ (depending also on the parameter
$\delta$). 

For $\lambda=0$ the solution of equation \eqref{eq:8} passing through
the origin is $y=0$, so $\sqrt{\Delta}|_{\lambda=0}=4-6x$. 
Similarly, for $x=0$ we have  $\sqrt{\Delta}|_{x=0}=4$.  
Therefore, by Hadamard Lemma we get 
$$\sqrt{\Delta}= 4-6x-\frac{27}{2} \delta \lambda x^2+ \lambda x\cdot o_{x,\lambda}(1) \ . $$

We deduce from \eqref{eq:8} the (smooth local) solution 
$$y=\dfrac{(4-6x)-\sqrt{\Delta}}{18\delta x} = \dfrac{3}{4} \lambda x+ \lambda \cdot o_{\lambda,x}(1) \ . $$
The graph of the corresponding $\lambda$-depending family of branches is parameterized by a function, whose $3$-jet is 
$$\left(x^2+\frac{3}{4} \lambda^2x,x^3,\lambda\right)\ .$$ 
This function is therefore $\A$-equivalent to the $\A$-sufficient $3$-jet $(x^2,x^3,\lambda)$, whose image is the semicubical cuspidal edge.
This ends the proof.
\end{proof}

The two $1$-parameter deformations $H$ and $K$ are obtained deforming a $(*)$-generic $1$-flat CTF of support $\Gamma$ among $(*)$-generic $1$-flat CTFs with the same support. 
Given such a family, there are indeed two such deformations, according to which cusp of the envelope we take as support. 

In the general theory of tangential families (see \cite{stable}, \cite{simple}), these deformations are called {\it tangential}.

\begin{theorem}
Consider a tangential deformation of $\psi_\delta$ (that is $H$ or $K$). 
For every small enough value of the deformation parameter, the two envelope branches have a second order tangency.
\end{theorem}

\begin{proof}
We consider the deformation $H$. 
The computations carried out in the proof of Theorem \ref{thm:5} show that the two branches of the critical set intersect themselves at the origin $x=0,y=0$.  
For $\lambda\not=0$ fixed, a suitable change of the coordinate $y$ brings the germ $\psi_\delta$ to the form 
$$\left(\lambda y+o_{x,y}(3),y^2-\frac{2}{\lambda}y^3-\frac{2}{\lambda}x^2y+x^3+o_{x,y}(3)\right)  
 \ .$$
Now, working modulo $\m_{x,y}^4\times\m_{x,y}^4$ we have 
$$\left(\lambda y ,y^2-\frac{2}{\lambda}(y^3-x^2y)+x^3\right)  \sim \left(y,x^3-\frac{2}{\lambda}x^2y\right)  \sim (y,x^3+x^2y)
 \ .$$
The latter $3$-jet is $\A$-sufficient (see \cite{rieger}); it is the normal form of the second order self tangency, which is the only stable local singularity of tangential families (see \cite{stable}).
\end{proof}

\begin{figure}[h]
  \centering
   \scalebox{.45}{\input{fig_3_flat.pstex_t}}
  \caption{Perestroikas of the double cusp singularity.}
 \label{fig:3}
\end{figure}
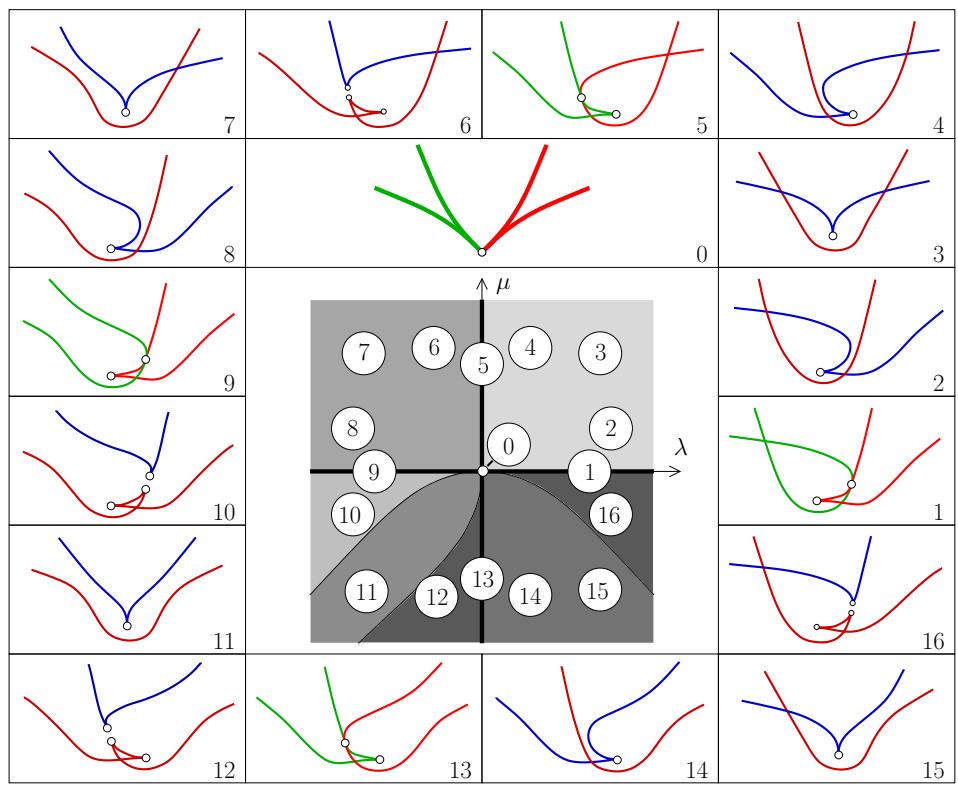

The second order self tangency of envelopes is stable under tangential deformations; under non tangential deformations the envelope experiences a beaks metamorphosis. 
The complete bifurcation diagram of the double cusp envelope
singularity is shown in figure \ref{fig:3} for $\delta>0$ (a change of
sign of $\delta$ inverses the orientation of the $\mu$ axis). 

\begin{remark*}
The bifurcation diagram of the double cusp singularity can be
interpreted also as the bifurcation diagram of the apparent contour of
a vertical Whitney Umbrella under small deformations of the direction
of projection. 
All the regular curves on the Umbrella surface, passing through the
singular point, have parallel velocities at this point. 
The ``verticality'' of the Whitney Umbrella means that the projection
fiber passing trough the Umbrella singularity is tangent to this
characteristic direction.  
\end{remark*}


\noindent{\sc Gianmarco Capitanio}\\ 
Department of Mathematical Sciences,\\
Mathematics and Oceanography Building,\\
Peach Street,\\
Liverpool L69 7ZL \\ 
United Kingdom

\end{document}

%% file: fig_1_flat.pstex_t
\begin{picture}(0,0)%
\includegraphics{fig_1_flat.pstex}%
\end{picture}%
\setlength{\unitlength}{3947sp}%
\begingroup\makeatletter\ifx\SetFigFont\undefined%
\gdef\SetFigFont#1#2#3#4#5{%
  \reset@font\fontsize{#1}{#2pt}%
  \fontfamily{#3}\fontseries{#4}\fontshape{#5}%
  \selectfont}%
\fi\endgroup%
\begin{picture}(5365,2186)(-14,-4035)
\put(4544,-2241){\makebox(0,0)[b]{\smash{{\SetFigFont{17}{20.4}{\rmdefault}{\mddefault}{\updefault}intersection}}}}
\put(4544,-2514){\makebox(0,0)[b]{\smash{{\SetFigFont{17}{20.4}{\rmdefault}{\mddefault}{\updefault}with the plane}}}}
\put(4551,-2807){\makebox(0,0)[b]{\smash{{\SetFigFont{17}{20.4}{\rmdefault}{\mddefault}{\updefault}$Z=0$}}}}
\end{picture}%

%% file: fig_2_flat.pstex_t
\begin{picture}(0,0)%
\includegraphics{fig_2_flat.pstex}%
\end{picture}%
\setlength{\unitlength}{3947sp}%
\begingroup\makeatletter\ifx\SetFigFont\undefined%
\gdef\SetFigFont#1#2#3#4#5{%
  \reset@font\fontsize{#1}{#2pt}%
  \fontfamily{#3}\fontseries{#4}\fontshape{#5}%
  \selectfont}%
\fi\endgroup%
\begin{picture}(20285,11557)(-12680,-2441)
\put(7510,6295){\makebox(0,0)[b]{\smash{{\SetFigFont{25}{30.0}{\rmdefault}{\mddefault}{\updefault}$\Gamma$}}}}
\put(7515,2327){\makebox(0,0)[b]{\smash{{\SetFigFont{25}{30.0}{\rmdefault}{\mddefault}{\updefault}$\Gamma$}}}}
\put(6451,4739){\makebox(0,0)[b]{\smash{{\SetFigFont{25}{30.0}{\rmdefault}{\mddefault}{\updefault}$\phi_\xi$}}}}
\put(7540,-1583){\makebox(0,0)[b]{\smash{{\SetFigFont{25}{30.0}{\rmdefault}{\mddefault}{\updefault}$\Gamma$}}}}
\put(5111,514){\makebox(0,0)[b]{\smash{{\SetFigFont{25}{30.0}{\rmdefault}{\mddefault}{\updefault}$\phi_\xi$}}}}
\put(6133,8816){\makebox(0,0)[b]{\smash{{\SetFigFont{25}{30.0}{\rmdefault}{\mddefault}{\updefault}$\phi_\xi$}}}}
\put(-11928,7934){\makebox(0,0)[b]{\smash{{\SetFigFont{25}{30.0}{\rmdefault}{\mddefault}{\updefault}$\phi_\xi$}}}}
\put(-11849,4064){\makebox(0,0)[b]{\smash{{\SetFigFont{25}{30.0}{\rmdefault}{\mddefault}{\updefault}$\phi_\xi$}}}}
\put(-7188,-399){\makebox(0,0)[b]{\smash{{\SetFigFont{25}{30.0}{\rmdefault}{\mddefault}{\updefault}$\phi_\xi$}}}}
\put(-6913,6323){\makebox(0,0)[b]{\smash{{\SetFigFont{25}{30.0}{\rmdefault}{\mddefault}{\updefault}$\Gamma$}}}}
\put(-6908,2295){\makebox(0,0)[b]{\smash{{\SetFigFont{25}{30.0}{\rmdefault}{\mddefault}{\updefault}$\Gamma$}}}}
\put(-6873,-1605){\makebox(0,0)[b]{\smash{{\SetFigFont{25}{30.0}{\rmdefault}{\mddefault}{\updefault}$\Gamma$}}}}
\put(-3524,-2311){\makebox(0,0)[b]{\smash{{\SetFigFont{29}{34.8}{\rmdefault}{\mddefault}{\updefault}$\epsilon=0$}}}}
\put(-3524,1589){\makebox(0,0)[b]{\smash{{\SetFigFont{29}{34.8}{\rmdefault}{\mddefault}{\updefault}$\epsilon=0$}}}}
\put(-3524,5489){\makebox(0,0)[b]{\smash{{\SetFigFont{29}{34.8}{\rmdefault}{\mddefault}{\updefault}$\epsilon=0$}}}}
\put(359,6280){\makebox(0,0)[b]{\smash{{\SetFigFont{25}{30.0}{\rmdefault}{\mddefault}{\updefault}$\Gamma$}}}}
\put(-2772,8111){\makebox(0,0)[b]{\smash{{\SetFigFont{25}{30.0}{\rmdefault}{\mddefault}{\updefault}$\phi_0$}}}}
\put(364,2312){\makebox(0,0)[b]{\smash{{\SetFigFont{25}{30.0}{\rmdefault}{\mddefault}{\updefault}$\Gamma$}}}}
\put(-2774,4139){\makebox(0,0)[b]{\smash{{\SetFigFont{25}{30.0}{\rmdefault}{\mddefault}{\updefault}$\phi_0$}}}}
\put(389,-1598){\makebox(0,0)[b]{\smash{{\SetFigFont{25}{30.0}{\rmdefault}{\mddefault}{\updefault}$\Gamma$}}}}
\put(-2768,236){\makebox(0,0)[b]{\smash{{\SetFigFont{25}{30.0}{\rmdefault}{\mddefault}{\updefault}$\phi_0$}}}}
\end{picture}%

%% file: fig_3_flat.pstex_t
\begin{picture}(0,0)%
\includegraphics{fig_3_flat.pstex}%
\end{picture}%
\setlength{\unitlength}{3947sp}%
\begingroup\makeatletter\ifx\SetFigFont\undefined%
\gdef\SetFigFont#1#2#3#4#5{%
  \reset@font\fontsize{#1}{#2pt}%
  \fontfamily{#3}\fontseries{#4}\fontshape{#5}%
  \selectfont}%
\fi\endgroup%
\begin{picture}(13281,10824)(-11,-9973)
\put(12976,-886){\makebox(0,0)[b]{\smash{{\SetFigFont{20}{24.0}{\rmdefault}{\mddefault}{\updefault}$4$}}}}
\put(4996,-7411){\makebox(0,0)[b]{\smash{{\SetFigFont{20}{24.0}{\rmdefault}{\mddefault}{\updefault}$11$}}}}
\put(5976,-7471){\makebox(0,0)[b]{\smash{{\SetFigFont{20}{24.0}{\rmdefault}{\mddefault}{\updefault}$12$}}}}
\put(6616,-7231){\makebox(0,0)[b]{\smash{{\SetFigFont{20}{24.0}{\rmdefault}{\mddefault}{\updefault}$13$}}}}
\put(12976,-6286){\makebox(0,0)[b]{\smash{{\SetFigFont{20}{24.0}{\rmdefault}{\mddefault}{\updefault}$1$}}}}
\put(12976,-4486){\makebox(0,0)[b]{\smash{{\SetFigFont{20}{24.0}{\rmdefault}{\mddefault}{\updefault}$2$}}}}
\put(12976,-2686){\makebox(0,0)[b]{\smash{{\SetFigFont{20}{24.0}{\rmdefault}{\mddefault}{\updefault}$3$}}}}
\put(12901,-8086){\makebox(0,0)[b]{\smash{{\SetFigFont{20}{24.0}{\rmdefault}{\mddefault}{\updefault}$16$}}}}
\put(12901,-9886){\makebox(0,0)[b]{\smash{{\SetFigFont{20}{24.0}{\rmdefault}{\mddefault}{\updefault}$15$}}}}
\put(9601,-9886){\makebox(0,0)[b]{\smash{{\SetFigFont{20}{24.0}{\rmdefault}{\mddefault}{\updefault}$14$}}}}
\put(6301,-9886){\makebox(0,0)[b]{\smash{{\SetFigFont{20}{24.0}{\rmdefault}{\mddefault}{\updefault}$13$}}}}
\put(3001,-9886){\makebox(0,0)[b]{\smash{{\SetFigFont{20}{24.0}{\rmdefault}{\mddefault}{\updefault}$12$}}}}
\put(3001,-8086){\makebox(0,0)[b]{\smash{{\SetFigFont{20}{24.0}{\rmdefault}{\mddefault}{\updefault}$11$}}}}
\put(6976,-5366){\makebox(0,0)[b]{\smash{{\SetFigFont{20}{24.0}{\rmdefault}{\mddefault}{\updefault}$0$}}}}
\put(8111,-5726){\makebox(0,0)[b]{\smash{{\SetFigFont{20}{24.0}{\rmdefault}{\mddefault}{\updefault}$1$}}}}
\put(8406,-5111){\makebox(0,0)[b]{\smash{{\SetFigFont{20}{24.0}{\rmdefault}{\mddefault}{\updefault}$2$}}}}
\put(8261,-4061){\makebox(0,0)[b]{\smash{{\SetFigFont{20}{24.0}{\rmdefault}{\mddefault}{\updefault}$3$}}}}
\put(7276,-3981){\makebox(0,0)[b]{\smash{{\SetFigFont{20}{24.0}{\rmdefault}{\mddefault}{\updefault}$4$}}}}
\put(6621,-4231){\makebox(0,0)[b]{\smash{{\SetFigFont{20}{24.0}{\rmdefault}{\mddefault}{\updefault}$5$}}}}
\put(5936,-3991){\makebox(0,0)[b]{\smash{{\SetFigFont{20}{24.0}{\rmdefault}{\mddefault}{\updefault}$6$}}}}
\put(4956,-4061){\makebox(0,0)[b]{\smash{{\SetFigFont{20}{24.0}{\rmdefault}{\mddefault}{\updefault}$7$}}}}
\put(4796,-5111){\makebox(0,0)[b]{\smash{{\SetFigFont{20}{24.0}{\rmdefault}{\mddefault}{\updefault}$8$}}}}
\put(5096,-5726){\makebox(0,0)[b]{\smash{{\SetFigFont{20}{24.0}{\rmdefault}{\mddefault}{\updefault}$9$}}}}
\put(4756,-6336){\makebox(0,0)[b]{\smash{{\SetFigFont{20}{24.0}{\rmdefault}{\mddefault}{\updefault}$10$}}}}
\put(7276,-7456){\makebox(0,0)[b]{\smash{{\SetFigFont{20}{24.0}{\rmdefault}{\mddefault}{\updefault}$14$}}}}
\put(8216,-7396){\makebox(0,0)[b]{\smash{{\SetFigFont{20}{24.0}{\rmdefault}{\mddefault}{\updefault}$15$}}}}
\put(8376,-6336){\makebox(0,0)[b]{\smash{{\SetFigFont{20}{24.0}{\rmdefault}{\mddefault}{\updefault}$16$}}}}
\put(9376,-5386){\makebox(0,0)[b]{\smash{{\SetFigFont{20}{24.0}{\rmdefault}{\mddefault}{\updefault}$\lambda$}}}}
\put(9676,-2686){\makebox(0,0)[b]{\smash{{\SetFigFont{20}{24.0}{\rmdefault}{\mddefault}{\updefault}$0$}}}}
\put(3076,-2686){\makebox(0,0)[b]{\smash{{\SetFigFont{20}{24.0}{\rmdefault}{\mddefault}{\updefault}$8$}}}}
\put(3076,-886){\makebox(0,0)[b]{\smash{{\SetFigFont{20}{24.0}{\rmdefault}{\mddefault}{\updefault}$7$}}}}
\put(3076,-4486){\makebox(0,0)[b]{\smash{{\SetFigFont{20}{24.0}{\rmdefault}{\mddefault}{\updefault}$9$}}}}
\put(3001,-6286){\makebox(0,0)[b]{\smash{{\SetFigFont{20}{24.0}{\rmdefault}{\mddefault}{\updefault}$10$}}}}
\put(6376,-886){\makebox(0,0)[b]{\smash{{\SetFigFont{20}{24.0}{\rmdefault}{\mddefault}{\updefault}$6$}}}}
\put(9676,-886){\makebox(0,0)[b]{\smash{{\SetFigFont{20}{24.0}{\rmdefault}{\mddefault}{\updefault}$5$}}}}
\put(6901,-3061){\makebox(0,0)[b]{\smash{{\SetFigFont{20}{24.0}{\rmdefault}{\mddefault}{\updefault}$\mu$}}}}
\end{picture}%